\documentclass[11pt,leqno]{article}
\usepackage{amsthm,amsfonts,amssymb,epsfig,graphics,amsmath,oldgerm}
\usepackage[latin1]{inputenc}\relax

\newcommand{\RR}{{\mathbb R}}
\newcommand{\NN}{{\mathbb N}}

\newcommand{\cD}{\mathcal{D}}

\newcommand{\cL}{\mathcal{L}}

\newcommand{\bE}{\mathbf{E}}

\newcommand{\bG}{\mathbf{G}}
\newcommand{\bK}{\mathbf{K}}

\newcommand{\eps}{\varepsilon }
\newcommand{\vp}{\varphi}
\newcommand{\D}{\partial }

\renewcommand{\div}{\,\mathrm{ div}\,}
\newcommand{\rot}{\, \mathrm{curl}\,}
\newcommand{\curl}{\, \mathrm{curl}\,}
\newcommand{\na}{{\nabla}}

\newcommand{\re}{\mathrm{Re\,}}
\newcommand{\tg}{\mathrm{tg}}

\newtheorem{theo}{Theorem}[section]
\newtheorem{prop}[theo]{Proposition}

\newtheorem{lem}[theo]{Lemma}
\newtheorem{defi}[theo]{Definition}

\newtheorem{rem}[theo]{Remark}

\numberwithin{equation}{section}

\title{Global existence and uniqueness for the lake equations
       with vanishing topography :
       elliptic estimates for degenerate equations}

\author{Didier {\sc Bresch} \footnote{LMC-IMAG, (CNRS-INPG-UJF)
38051 Grenoble cedex, France. Email: didier.bresch@imag.fr}
and  Guy {\sc M\'etivier} \footnote{MAB et IUF,  Universit\'e de Bordeaux I,
33405 Talence cedex, France. Email : Guy.Metivier@math.u-bordeaux.fr}}

\begin{document}

\maketitle


\begin{abstract}

\smallskip

  This paper deals with global existence and uniqueness for the lake equations with
a bottom topography vanishing on the shore.
  Our result generalizes previous studies that assumed the depth to be nondegenerate.
  Elliptic estimates for degenerate equations are established studying the behavior
of the associated Green function.

\bigskip



\end{abstract}


\section{Introduction}
  In this paper we are interested in a two-dimensional geophysical model that has been
essentially described by Greenspan in \cite{Gr} page 235.
  This system describes the evolution of the vertically averaged horizontal component $v(t,x,y)$
of a three dimensional velocity vector $u(t,x,y,z)$ and reads
\begin{equation}
\label{Lake}
\left\{
\begin{aligned}
& \partial_t (bv) + {\rm div}_{\bf x} (b v \otimes v) + b\nabla_{\bf x} p = 0 \hbox{ in } \Omega, \\
& {\rm div}_{\bf x}(bv)= 0 \hbox{ in } \Omega, \\
& (bv)\cdot n = 0 \hbox{ on } \partial \Omega, \\
& (bv)|_{t=0} = m_0 \hbox{ in } \Omega
\end{aligned}\right.
\end{equation}
where ${\bf x}$ denotes the horizontal components ${\bf x} = (x,y)$.
  The lake equations may be seen as the low Froude number limit of the usual inviscid
shallow water equations when the initial height converges to a nonconstant function
depending on the space variable, namely $b({\bf x})$, see for instance \cite{BrGiLi}.
  Recall that the inviscid shallow water reads
\begin{equation}
\label{SW}
\left\{
\begin{aligned}
& \partial_t h + {\rm div}_{\bf x}(hv) = 0 \hbox{ in } \Omega, \\
& \partial_t (hv) + {\rm div}_{\bf x} (h v \otimes v)
                  + h\nabla_{\bf x}\frac{(h-b)}{\rm Fr^2} = 0 \hbox{ in } \Omega, \\
& (hv)\cdot n = 0 \hbox{ on } \partial \Omega, \\
& (hv)|_{t=0} = m_0, \qquad h|_{t=0} = h_0 \hbox{ in } \Omega,
\end{aligned}\right.
\end{equation}
where ${\rm Fr}$ is the Froude number, $h$ is total depth, $v$ the vertical mean
value of the horizontal components of the velocity.
\medskip

\noindent {\it The constant case $b= {\rm const}$: The 2D incompressible Euler equations.}
  We remark that the case $b\approx 1$ reduces to the standart 2D incompressible Euler system.
  The Cauchy problem for incompressible Euler equations is very well understood
and we refer the reader to various existing surveys on the question : see for instance
\cite{Yo}, \cite{Ma}, \cite{Ch} and \cite{MaBe}.
  Concerning the bounded domain case, we refer to \cite{Li} where the method used
by Youdovitch is explained and various results described.
  The approach is the construction of the solutions as the inviscid limit of
solutions to a Navier-Stokes system with artificial viscosity and boundary conditions
$v\cdot n = 0$ and $\omega = 0$ where $\omega = \partial_y v_1 - \partial_x v_2$ is the vorticity.
   We also mentionned \cite{ClMiRo} where a wall law related to the cauchy stress tensor replaces
the null condition on the vorticity.
   In all these works, to get the result, elliptic estimates are used on the following equation
$$ -\Delta_{\bf x} \Psi = \omega \hbox{ in } \Omega, \qquad
    \Psi|_{\partial \Omega} = 0.
$$
 Namely the following estimate
$$ \|\Psi\|_{W^{2,p}(\Omega)} \le C p \,\|\omega\|_{L^p(\Omega)}
$$
where $C$ does not depend on $p$.
 Such elliptic estimates with nondegenerate coefficients follow
from \cite{ADN1}--\cite{ADN2}.
 We also mention in $\RR^2$ the method based on approximate solutions which
strongly uses the following property on $\nabla_{\bf x} v$, see for instance \cite{MaBe}:
Let $v$ be a smooth, $L^2(\RR^2)\cap L^\infty(\RR^2)$, divergence free field
and let $\omega = {\rm curl} v$, then
$$\|\nabla v\|_{L^\infty(\RR^2)} \le c (1+ \ln^+ \|v\|_{L^3(\RR^2)}
               + \ln^+ \|\omega\|_{L^2(\RR^2)})(1+\|\omega\|_{L^\infty(\RR^2)})$$
where $\ln^+ g = 0$ if $g\le 1$ and $\ln^+ g = \ln g$ if $g > 1$.

\medskip

\noindent {\it The nonconstant case $b \not = {\rm const}$: The lake equations.}
  If the depth is not assumed to be constant, we have to define a new vorticity
variable to get a transport equation on this quantity.
  Namely defining
$$ \omega = \frac{1}{b}(\partial_y v_1 - \partial_x v_2)
$$
we get that it is transported along the flow $v$ namely
\begin{equation}
\label{transp}
 \partial_t \omega + v\cdot \nabla_{\bf x} \omega = 0.
 \end{equation}
Together with the incompressibility condition $\div (b v) = 0$, this implies, formally,  conservation of  the $L^p ( b dx)$ norms of $\omega$, for $p \in [1, \infty]$.

  The case $b\ge C > 0$ has been studied in \cite{LOT} proving global
existence and uniqueness for an initial vorticity $\omega_0$ in $L^\infty(\Omega)$
following Youdovitch's method.
  The regularity of the lake equations has been also addressed in \cite{Ca} in this
case but where $|\nabla b(x)|\to \infty$ as $x$ tends to the boundary.
 Here we study the degenerate case, namely the case when $b$ stricly positive
in $\Omega$ vanishing on the shore $\partial \Omega$.
 Note that, even if the model seems to be not physically relevant in the presence of beaches,
similar models are used by several geophysicits.
 In a forthcoming paper, the authors will try to better understood the right
model to be considered in the presence of beaches.

 We will assume that $b > 0$ on $\Omega$ and  $b = 0$ on $\D \Omega$.
 More precisely, we assume that $\Omega$ is bounded with smooth boundary
$\D \Omega$, and
  \begin{equation}
\label{b}
b  = \vp^a
\end{equation}
where $ a > 0$ and
$\Omega = \{ \vp > 0 \}$  with $\vp \in C^\infty(\overline \Omega)$ and
$\na_{\bf x} \vp \ne 0$ on $\D \Omega $.

  In order to get our result, an elliptic regularity result is required for a degenerate equation on the associated stream-function.
  This is the main part of the paper.
It is strongly related to a careful study of the associated Green function given
in \cite{GS}--\cite{GS1}.

 The paper is organized as follows. The main results will be summarized in Section 2.
 In section 3 we will mention the existence and uniqueness result and the key result of the paper: Namely the elliptic estimates on the stream function $\Psi$.
 Sections 4, 5, 6 are respectively devoted to the
proof of such $L^p$ estimates by proving the normal regularity and vanishing
property, by showing the H\"older estimates and finally by establishing the $L^p$ estimates respectively.
  The well posedness result will be proved in the last section adapting Youdovitch's proof used for the standart 2D incompressible Euler equations.

  For the sake of simplicity, throughout the paper, we will supress the ${\bf x}$ index
from the partial derivatives.


\section{Main results.}
   The goal of the paper is to prove an existence and uniqueness result of a global
weak solution for the inviscid equations \ref{Lake}. We assume that 
$\Omega$ is  a smooth simply connected bounded domain and that 
$b$ is given by \eqref{b}.
 As usual, we eliminate the pressure
$p$ from the equation and consider the 
\emph{vorticity-stream formulation of the 2D lake equation}.  The weak form of the transport equation \eqref{transp} is
\begin{equation}
 \label{transpf} \D_t (b \omega)  + \div ( b \omega v)  = 0.
 \end{equation}
 The vorticity and and the velocity are linked by the equations: 
 \begin{equation}
\label{vort}
\div (b v ) = 0, \quad  \curl v = b \omega, \quad  ( b v) \cdot n _{| \D \Omega}  = 0. 
\end{equation}
This leads, see for instance \cite{MaPu}, to introduce a stream function $\Psi$ such that
\begin{equation}
\label{stram}
v  = \frac{1}{b} \na^\perp \Psi =  \frac{1}{b}   (\D_2 \Psi , - \D_1 \Psi),
\qquad \Psi|_{\partial\Omega} = 0.
\end{equation}
 In the non-simply connected case, which is discussed later in a 
remark, one has to specify different boundary values on each 
connected component of $\partial\Omega$, see \cite{MaPu}.
 The model for $\Psi$ then reads
\begin{equation}
\label{streameq}
\div ( \frac{1}{b} \na \Psi)  = b \omega, \quad  \Psi_{| \D \Omega}  = 0. 
\end{equation}
Note that, by Hardy's inequalities, the space 
$C^\infty_0 (\Omega)$ is dense in the   space 
$\{ \Psi \in H^1_0(\Omega) ; b^{ - 1/2} \na \Psi \in L^2(\Omega) \}$. Therefore, for 
$f \in L^2(\Omega)$ the problem  
 \begin{equation}
\label{eqpsipsi}
\div ( \frac{1}{b} \na \Psi)  = f, \quad  \Psi_{| \D \Omega}  = 0, 
\end{equation}
has a unique solution $\Psi $ in this space, denoted by 
$\Psi = K f$. 
   We make the following definition.
\begin{defi} 
\label{def21}
Given $\omega_0\in L^\infty(\Omega)$, $(v,\omega)$ is a weak solution to the
vorticity-stream formulation of the 2D lake equation with initial data $\omega_0$, provided

 {\rm i)}      $\omega  \in L^\infty ([0,T])\times \Omega)$ and 
 $b\omega \in C^0([0,T]; L^\infty_{w*}(\Omega))$,
 
  {\rm ii)} $ b v =  \na^\perp K (b \omega) \in C^0([0, T]; L^2(\Omega)  $,
  
  {\rm iii)} For all $\varphi \in C^1([0,T]\times \overline \Omega)$ and
$t_1\in [0,T]$:
$$ \int_\Omega b\varphi(t_1,x)\omega (t_1,x) \, dx
  - \int_\Omega b\varphi(0,x) \omega_0(x) \, dx = \int_{0}^{t_1}
  \int_\Omega (b\partial_t \varphi + b v\cdot \nabla \varphi) \omega \, dx dt.
$$

\end{defi}

\begin{theo}
\label{th22}
\noindent {\rm i)}    {\bf (Regularity) }  Assume that $\omega\in L^\infty((0,T)\times\Omega)$ and
$  b  v\in L^\infty([0,T] ; L^2(\Omega))$ such that ${\rm div}(bv)=0$ satisfy the weak formulation.
 Then $\omega \in C^0([0,T]; L^r(\Omega))$ and $v\in C^0([0,T]; W^{1,r}(\Omega))$ for
all $r<+\infty$. 
Moreover, there is $C$ such that for all $p\ge 3 $,
$$
 \|\nabla v\|_{L^p(\Omega)} \le  C  p \|b\omega\|_{L^p (\Omega)}.
 $$
 In addition, we get the following boundary condition on $v$  
$$
 v\cdot n = 0 \hbox{ on } \partial \Omega.
$$

  {\rm ii)} {\bf (Existence)} For all  $\omega_0 \in L^\infty(\Omega)$, 
  there exists a global weak solution $(v, \omega)$ to the
vorticity-stream formulation of $\eqref{Lake}$.

  {\rm iii)} {\bf (Uniqueness)} The weak   vorticity-stream solution is  unique. 
\end{theo}

 This result follows the Youdovitch's procedure in constructing the solution
as the inviscid limit of solutions of a system with artificial viscosity
which is the analog of Navier--Stokes with respect to the Euler equations.
 The core of the proof is a regularity result for a system degenerate due to
coefficient vanishing on the boundary of the domain.
 More precisely, the main part of the paper will concern the following
main result.
 Let $\Omega$ be a smooth simply connected bounded domain in $\RR^n$ with $n\ge 1$
and $b$ defined by \eqref{b}.
 Consider the system:
\begin{equation}
\label{16}
\left\{
\begin{aligned}
 & \div (b v)  = 0 \hbox{ in } \Omega,  \quad  (bv) \cdot n _{ | \D \Omega} = 0   ,  \\
 & \curl v =  f \hbox{ in } \Omega
\end{aligned}\right.
\end{equation}
for
\begin{equation}
\label{17}
\left\{ \begin{aligned}
 &   b v \in L^2(\Omega),  \\
&    f \in L^p (\Omega).
\end{aligned}\right.
\end{equation}

 Using the definition of $b$ given in \eqref{b}, we prove the following result on which the existence and uniqueness result is based.

\begin{theo}
\label{mainest}
If $(v,f)$ satisfy $\eqref{16}$ and $ \eqref{17}$ with $p \in ]2 , \infty[$,  then
\begin{equation}
\label{19}
  v \in   C^{1 - n/p }(\overline \Omega), \quad
  \na v\in  L^p(\Omega) .
\end{equation}
There is a constant $C_p  $ independent of $(u, \omega)$ such that
\begin{equation}
\label{110a}
\Vert  v \Vert_{C^{1 - n/p }(\overline \Omega)}    \le C_p
\big( \Vert   f  \Vert_{L^p}  + \Vert  b v \Vert_{L^2} \big).
\end{equation}
In addition 
\begin{equation}
\label{110b}
v \cdot n  = 0 \quad   on \ \D \Omega.
\end{equation}
Moreover, for all $p_0 > 2$,  there is a constant $C$ independent of $(v, f )$ and 
$p \in [p_0, \infty[ $  such that
\begin{equation}
\label{110}
  \frac{1}{p} \Vert \na v \Vert_{L^p }  \le C
\big( \Vert  f \Vert_{L^p}  + \Vert  b v\Vert_{L^2} \big).
\end{equation}
\end{theo}

In this statement, for $\mu \in ]0, 1[$,  $C^\mu (\overline \Omega)$ is the usual
space of continuous functions on $\overline \Omega$ which satisfy the H\"older condition
of order $\mu$.
  In particular, \eqref{110a} implies that
\begin{equation}
\label{110aa}
\Vert  v \Vert_{L^\infty( \Omega)}    \le C
\big( \Vert   f \Vert_{L^\infty}  + \Vert  b v \Vert_{L^2} \big).
\end{equation}

\begin{rem}
\textup{ It is remarkable that the final $W^{ 1, p}$ smoothness of $v$ is independent 
of $b$ in \eqref{16}.  When $p = 2$,   it is still true that $ v \in W^{1, 2}$ and indeed estimate
\eqref{110} is uniform for $p \in [2, \infty [$. This is shown using estimates of \cite{BG, BC}
in the proof below in place of those of  \cite{BCM}. However, in this paper we are mainly interested 
in these estimates as $p \to \infty$, and we leave aside  the case $p = 2$.  Note also that \eqref{110a} follows from \eqref{110} and the 
Sobolev embedding $W^{1, p} \subset C^\mu$ for $\mu = 1 - 2/p > 0$. 
The proof of the uniform estimate  \eqref{110} given below relies on proving \eqref{110a} for a given $p_0 > n$ first, and on winning the final $L^p$ estimate
for $\na v$, then. This is why we have stated them independently.  }
\end{rem}

\begin{rem}
\textup{ We  include here some  remarks on the link between the vorticity-stream formulation 
and  \eqref{Lake}.}

\textup{1.  Suppose first  that   $v \in \mathcal{C}^0 ([0, T], H^1(\Omega))$ is a solution  in the sense  of distributions of \eqref{Lake}. For $u \in H^1$, the following identity holds:   }
$$
\curl \big( \frac{1}{b} \div ( b u \otimes u) \big)  =  \div \big( (\curl u) u \big)  -
\curl \big( \frac{1}{b} ( \div  b u)   u  \big) 
$$
\textup{implying that $\omega = \curl v$ satisfies transport equation \eqref{transpf} in the sense of distributions thus condition $iii)$ in Definition~$\ref{def21}$.  Moreover, \eqref{vort} is satisfied by
\eqref{Lake} and by definition of $\omega$. Thus $v = b^{-1} \na^\perp \Psi$, (see the beginning of section 3) and the stream function $ \Psi \in H^{1}_0 $ satisfies \eqref{streameq}, 
hence $v = b^{-1} \na^\perp K (b \omega)$ as in $ii)$ of Definition~$\ref{def21}$.  }

\textup{2.  Conversely, suppose that $(\omega, v)$   is a weak solution of the
vorticity-stream formulation. By part $i)$ of Theorem~$\ref{th22}$ 
$v \in \mathcal{C}^0 ([0, T], L^2(\Omega))$ and, because $\omega$ satisfies
\eqref{transpf}, we also have $\D_t \omega \in L^\infty([0, T], H^{-1}(\Omega))$. In particular, }
$$
\div (b \D_t  v ) = 0, \quad  \curl (\D_t v)  = b \D_t \omega \in   L^\infty([0, T], H^{-1}(\Omega)), 
$$
\textup{implying that $\D_t v  \in  L^\infty([0, T], L^2_{loc}(\Omega))$. The same computation 
as in 1), implies that in the sense of distributions }
$$
\curl  \Big( \D_t  v+ \frac{1}{b}\div (b v \otimes v) \big) \Big)
= \D_t(b\omega) + \div ( b v \omega) = 0. 
$$
\textup{Therefore, assuming that $\Omega$ is simply connected, there is a pressure
$p \in L^\infty([0, T], \mathcal{D}' (\Omega))$ such that  \eqref{Lake} is satisfied.} 

\end{rem}

\begin{rem} 
\textup{
When $\Omega$ has some islands that means it is a 
non-simply connected domain, some generalized circulations are
specified in order to uniquely determine the velocity $v$, see \cite{MaPu}.
 More precisely let $\Omega$ with $n$ islands $I_1$, $\dots$, $I_n$, the missing parameters are the boundary values $\lambda_i$ with $i=1,\dots, n$ of the stream 
function on each of the boundary component $\D I_1$, $\dots$, $\D I_n$.
 Assuming the generalized circulations to be zero, the velocity is
given by $v = v_0 + \sum_{i=1}^n \lambda_i v_i$ with
$v_j = (\nabla^\bot \Psi_j) / b$ where $\Psi_j$ ($j=0,\dots,n$) 
are the unique solution of} 
$$ -{\rm div}(\frac{1}{b}\nabla \Psi_i) = 0, \qquad 
\Psi_i|_{\D I_j} = \delta_{i,j} \quad j=0,\dots,n \quad i=1,\dots,n.$$
\textup{and}
$$-{\rm div}(\frac{1}{b} \nabla \Psi_0) = b\omega, 
\qquad \Psi_0|_{\D \Omega} = 0.$$
\textup{  The $n$ coefficients $\lambda_j$ are uniquely determined using the
zero generalized circulations conditions.
  This decomposition shows that our study concerning the simply
connected domain may be applied to the non-simply connected one, 
the $L^p$ estimate remaining true.}  
\end{rem}

\section{Localization}
  Since $b v\in  L^2(\Omega)$ satisfies $\div (bv) = 0$
and $ (bv) \cdot n = 0  $ on $\D \Omega$,
there is a unique potential $\Psi $ such that
\begin{equation}
\label{upsi}
v   = \frac{1}{b} \na^\perp \Psi =
\frac{1}{b} ( \D_y \Psi, - \D_x \Psi) , \quad     \Psi_{|  \D \Omega} = 0. 
\end{equation}
Indeed, $\Psi$ is determined 
in $H^1_0(\Omega)$ as the solution of 
\begin{equation}
\label{eqpot1}
- \Delta \Psi = \rot ( b v  ) , \quad  \Psi_{|  \D \Omega} = 0.
\end{equation}
 There holds 
$$
\rot (b v) = b \rot v + \na b \times v= b f  - \frac{1}{b} \na b \cdot \na \Psi
$$
  With \eqref{b}, we obtain the equation
\begin{equation}
\label{eqpot2}
- \vp \Delta \Psi +  a \na \vp \cdot \na \Psi =  \vp^{ a+1} f  ,   
 \quad  \Psi_{|  \D \Omega} = 0.
\end{equation}
 
 At this stage, the dimension $d = 2$ plays no particular rule and we can
consider this equation on any smooth bounded domain $\Omega  = \{ \vp > 0 \}
\subset \RR^n$, with $\vp \in C^\infty(\RR^n )$ and $d \vp \ne 0$ on $\D \Omega $.
 One introduces the H\"older spaces
$C^\mu(\overline \Omega)$ and the Sobolev spaces 
$W^{k, p}(\Omega)$ based on $L^p$,   together with the weighted spaces 
$C_1^\mu(\overline \Omega)$ and  
$W_1^{k, p}(\Omega)$ of functions 
$u$ in $C^\mu(\overline \Omega)$ or   
$W^{k, p}(\Omega)$ such that $ \vp u$ in $C^{\mu+1} (\overline \Omega)$ or   
$W^{k+1, p}(\Omega)$ respectively. 

Writing 
\begin{equation*}
 \Psi =  \vp^{a+1} \Phi, \quad  u = \vp^{-a} \na^\perp \Psi = 
 \vp \na^\perp \Phi +  (a+1) \Phi  \na^\perp  \vp , 
 \end{equation*} 
Theorem~\ref{mainest}  follows from 

\begin{theo}
\label{mainest2}
If $\Psi \in H^1_0(\Omega)$ and $f   \in L^p(\Omega) $ with $p > n$ satisfy $\eqref{eqpot2}$, then
$   \Phi =  \vp^{-(a+1)} \Psi $ satisfies  
\begin{equation}
\label{mreg}
\Phi    \in \bigcap_{ \mu \le  1-n/p} C_1^{\mu} (\overline \Omega),
\quad   \Phi   \in \bigcap_{q\le p} W^{1, p}_1 (\Omega).
\end{equation}
Moreover, there are constants $C_\mu$, $\mu \in [0, 1-n/p]$ such that for all
$\Psi \in H^1_0 $ and $f   \in L^p  $ satisfying $\eqref{eqpot2}$,
\begin{equation}
\label{mest1}
  \Vert   \Phi  \Vert_{C^{\mu}_1 (\overline \Omega)}      \le  C_ \mu  
\big( \Vert  f  \Vert_{L^p}  + \Vert  \Psi \Vert_{H^1} \big),  
  \end{equation}
  For all $p_0 > n$, there is a constant $C_1$, such that  for all
$\Psi \in H^1_0 $ and $f   \in L^p  $ with $p \ge p_0$ satisfying $\eqref{eqpot2}$
  \begin{equation}
  \label{mest2}
    \frac{1}{p} \Vert \Phi \Vert_{W^{1, p}_1 }    \le  C_1 
\big( \Vert f  \Vert_{L^p}  + \Vert  \Psi \Vert_{H^1} \big).
\end{equation}

\end{theo}
  
 Equation \eqref{eqpot2} is a particular case of degenerate elliptic equations studied in
\cite{BG, BC, GS, BCM} and the estimates \eqref{mest1} \eqref{mest2} are mainly consequences of the
results obtained in these papers.
 Because it may be useful in other circumstances and also because it helps to understand the analysis,
 we will prove the estimates for a slightly  more general class of equations, namely
equations of the form
 \begin{equation}
\label{eqpot3}
L(x, \D_x) \Psi :=   -  \vp(x)     P_2 (x, \D_x ) \Psi     +   a  P_1(x, \D_x) \Psi  =  \vp^{a+1}  f 
\end{equation} 
where $a$ is a positive  real \emph{constant} and 
\begin{equation}
\label{56a}
P_2 = \sum_{j, k= 1}^n p_{j,k} (x) \D_j \D_k , \quad  P_1 = \sum_{j=1}^n p_{j} (x)  \D_j  
\end{equation}
have real and smooth coefficients on $\overline \Omega$, with $p_{j, k}(x)  = p_{k, j}(x)$.
 We further assume that $P_2$ is uniformly elliptic on $\overline \Omega$
\begin{equation}
\label{ellip}
\sum p_{j, k} (x) \xi_i \xi_k  > 0 \quad  \mathrm{for}  \ 
x \in \overline \Omega, \ \xi \in \RR^n\backslash\{0\}, 
\end{equation}and  
\begin{equation}
\label{Z}
  P_1(x, \xi)   =       \sum p_{j, k}(x)  \xi_j \D_k \vp (x)  \quad     \mathrm{on } \ \D \Omega    . 
\end{equation} 
 Notice that \eqref{ellip} and \eqref{Z} are independent of the choice of the defining function $\vp$
for $\Omega$.
 
 \begin{theo}
 \label{theomainest3}
 With assumptions as above, if $\Psi\in H^1_0(\Omega)$ and $f\in L^p(\Omega)$
satisfy $\eqref{eqpot3}$, then $\Phi = \vp^{-(a+1)} \Psi$ satisfies $\eqref{mreg}$ and there are
constants $C_\mu$ such that the  estimates $\eqref{mest1}$ and $\eqref{mest2}$
hold. 
 \end{theo}
 
\bigbreak
\noindent \emph{Reduction to a neighborhood of the boundary. }
  Standard elliptic theory implies  that on any open subset  $\Omega_1 \subset\!\!\subset \Omega$,
$\Psi$ belongs to   $W^{2, p} (\Omega_1)$ for all $p$ finite and thus to
$C^{\mu +1}( \overline  \Omega_1)$ for all $\mu \le  1-n/p$.
  Moreover,
\begin{equation}
\label{26}
\Vert    \Phi   \Vert_{C^{\mu +1}(\Omega_1)}   \le C_\mu \big(  
 \Vert f  \Vert_{L^p } + \Vert \Psi \Vert_{H^1}\big)
\end{equation}
and, for $p \in [2, \infty[$: 
\begin{equation}
\label{26b}
\frac{1}{p}   \Vert  \na^2  \Phi \Vert_{L^p(\Omega_1)}  \le C_1\big(  
 \Vert f  \Vert_{L^p } + \Vert \Psi \Vert_{H^1}\big)
\end{equation}
 Therefore, if $\chi \in C^\infty (\overline \Omega)$ is such that
$\chi - 1 \in C^\infty_0(\Omega)$, then
\begin{equation}
\label{eqpotloc}
 L (\chi \Psi ) = \vp^{a+1} \chi f  +  g  \quad  \Psi_{|  \D \Omega} = 0, 
\end{equation}
 where the commutator $g$ is $C^{1- \eps}$  by \eqref{26} and vanishes on a neighborhood 
of $\D \Omega$.
  Therefore, one can factor out any power of 
$\vp$, writing  $g = \vp^{a +1}  f_1$ with $f_1 \in L^\infty$ with norm bounded by the right hand
side of \eqref{26}.
  Therefore, it is sufficient to prove the estimates \eqref{mest1}--\eqref{mest2} for functions $\Phi$ that are supported in an arbitrary small neighborhood of $\D \Omega$.

\bigbreak
\noindent \emph{Local coordinates near the boundary.}
  The boundary $\D \Omega$ is a closed smooth manifold.
  Consider a coordinate patch $x' \mapsto \gamma(x')$ from an open set
$\omega \subset \RR^{n-1}$ to $\D \Omega$, with $\gamma \in C^\infty$ on $\overline \omega$.
   Taking $\nu (x')$ to be the inward normal to $\D \Omega$ at $\gamma (x' ) $,
we parametrize a neighborhood of $\gamma (\omega)$  by
$(x', x_n)  \in  \omega \times ]- \delta, \delta [ $ by considering the mapping
\begin{equation}
\label{loccoord}
\Gamma :  \  (x', x_n)   \mapsto   \gamma(x')  + x_n \nu(x'). 
\end{equation}
  In these coordinates, the equation \eqref{eqpot3} is transformed to
\begin{equation}
\label{eqloc}
 \widetilde L \Psi  = - x_n    \widetilde P_2 \Psi +   a  \widetilde P_1\Psi
                    =   x_n^{ a+1} f.
\end{equation}
   The ellipticity property \eqref{ellip} is preserved as well as \eqref{Z} which now reads
\begin{equation}
\label{Zloc}
  \widetilde P_1(x', 0 , \xi)  = \sum_{j=1}^n \tilde p_{j, n}(x', 0)  \xi_j.
\end{equation} 

 \bigbreak
 Collecting all these remarks, we see that Theorem~\ref{theomainest3} follows from the next 
estimates.
\begin{theo}
\label{mainestr} 
 If $\Psi \in H^1(\omega\times ]0, \delta [ )  $ and $f  \in L^p(\omega \times [0, \delta]) $ 
 with $p > n$ satisfy
 $\eqref{eqpotloc}$ and 
 \begin{equation}
\label{bcloc}
   \Psi_{| x = 0} = 0,  \quad  
\Psi _{|  \{ x > \delta /2\} }  = 0 . 
\end{equation} 
then 
$  \Psi =  x_n^{a+1} \Phi $ with 
\begin{equation}
\label{mregloc}
\Phi     \in    C_1^{1 - n/p} (\overline  \omega_1 \times [0, \delta] ),
\quad   \Phi   \in   W^{1, p}_1 (\omega_1 \times [0, \delta] )
\end{equation}
for all relatively compact open subset $\omega_1 \subset \omega$. 
  Moreover, there are constants $C_\mu$, $\mu \in [0, 1-n/p]$ such that for all
$\Psi \in H^1_0 $ and $f\in L^p $,  $p \ge p_0 > n$  satisfying $\eqref{eqpot2}$,
\begin{eqnarray}
\label{mest1loc}
  \Vert   \Phi  \Vert_{C^{\mu}_1 (\overline \omega_1 \times [0, \delta[ )}   &  \le &C_ \mu  
\big( \Vert  f  \Vert_{L^p}  + \Vert  \Psi \Vert_{H^1} \big),
  \\
  \label{mest2loc}
  \frac{1}{p} \Vert \Phi \Vert_{W^{1, p}_1(\omega_1 \times [0, \delta]) }  & \le& C_1 
\big( \Vert f  \Vert_{L^p}  + \Vert  \Psi \Vert_{H^1} \big).
\end{eqnarray}
\end{theo}

  Here $C^\mu_1$ [resp.  $W^{1, p}_1$] denote the spaces for functions
$u \in C^\mu$ [resp. $u \in W^{1, p}$] such that $ x_n u \in C^{\mu+1}$
[resp. $x_n u \in W^{2, p}$].


\section{Smoothness in normal direction and vanishing at the boundary }
  The first step in the proof of Theorem~\ref{mainestr} is to factor out
$x_n^{a+1}$ in $\Psi$ and to show that $\Phi $ is H\"older continuous in the normal 
variable.
  For simplicity, we drop the tildes in equation \eqref{eqpotloc}.

\begin{prop}
\label{prop32}
  With assumptions as in Theorem~$\ref{mainestr}$ there are $\mu > 0$, $C$
and a Banach space $E \subset \cD'(\omega)$, such that the function $\Phi = x_n^{-(a+1) } \Psi$  satisfies
on $\omega \times [0, \delta]$:
\begin{equation}
 \label{36}
 \Phi \in   C^{\mu } ([0, \delta]; E), \quad
 x_n  \D_{x_n} \Phi \in   C^{\mu } ([0, \delta] ; E).
\end{equation}
together with the estimates 
\begin{equation}
\label{estnorm}
\Vert   \Phi \Vert_{ C^\mu([0, \delta]; E) }  + 
\Vert x_n \D_{x_n}  \Phi \Vert_{ C^\mu([0, \delta]; E) }  \le  C \Big(
\Vert f  \Vert_{L^p}  + \Vert \Psi \Vert_{H^1}\Big) .
\end{equation}
\end{prop}

  Multiplying the equation \eqref{eqpotloc} by an appropriate factor, one can always assume that the coefficient of $\D_{x_n}^2$ is exactly $x_n$.
  Then, \eqref{Zloc} implies that the coefficient of $\D_{x_n}$ in $P_1$  is 1 when $x_n = 0$.
Thus \eqref{eqpotloc} reads
\begin{equation}
\label{31}
-  x_n \D_{x_n}^2  \Psi  +  a \D_{x_n} \Psi  =   x_n^{a+1} f  +  F , \quad   F = L' \Psi   , 
\end{equation}
with  
\begin{equation}
\label{erF}
L' :=   \sum_{j, k < n}   x_n   p_{j, k} \D_j \D_k   +   x_n  \sum_{j < n} p_{j, n} \D_j \D_{x_n} 
- a \sum_{j < n} p_j \D_y  -  a \tilde p_n x_n \D_{x_n} 
\end{equation}
where the   $\tilde p_n $ is smooth. 

  We consider functions of $x_n$ valued in (Banach) spaces of distributions in $x'$ :
$L^2([0, \delta]; E) $ , $C^\mu ([0, \delta], E)$  etc , with $E \subset \cD' (\omega)$.
The space $E$ is subject to change from line to line, and we note 
$L^2( \cD'_{\tg} ) $ , $C^\mu (  \cD'_{\tg})$  when the space is unspecified. 
   We also denote by  $x_n^\alpha \cD'_{\tg} $ the space of   products of
the function $x_n^ \alpha$ with arbitrary distributions in $\cD'(\omega)$.

  We know that
\begin{equation}
\label{32}
\Psi \in H^1 \subset  C^{1/2} (L^2) . 
\end{equation}
Using the equation \eqref{31}, we also see that 
$x_n \D_{x_n}^2 \Psi \in L^2(H^{-1})$ and hence 
\begin{equation}
\label{33}
x_n\D_{x_n} \Psi \in H^1(H^{-1} ) \subset   C^{1/2} (H^{-1}) .
\end{equation}

  We will make use of the next (classical) results (See also Proposition~4.1.2 in \cite{BCM} for
more general results).

\begin{lem}
\label{lembase}
For  $ \alpha > 0 $ and $\mu \in ]0, 1[$, the operators
\begin{equation}
\label{34}
I_\alpha u (x_n) =   \int_0^{x_n} \left( \frac{t}{x_n}\right) ^\alpha u(t) \frac{dt}{t} , \quad 
J_\alpha u (x_n) =      \int_{x_n}^\delta   \left( \frac{x_n}{t}\right) ^\alpha  u(t) \frac{dt}{t} 
\end{equation}
map continuously $C^\mu([0, \delta]; E)$ to $C^\mu([0, \delta]; E)$.
\end{lem}

\begin{proof}[Proof of Proposition~$\ref{prop32}$]
  Because $\Psi = 0$ for $x \ge \delta /2$, Equation \eqref{31} implies that
\begin{equation}
\label{37}
\D_x \Psi  = J_a ( x_n^{a+1} f +  F ) =   x_n^{a} \theta  +  J_a( F)  
\end{equation}
 with 
 \begin{equation}
\label{38}
\theta  = - \int_{x_n}^\delta f(t) dt  \in C^{1/2} (L^p).
\end{equation}
Because $\Psi \in H^1$ vanishes at $x_n= 0$, there holds
\begin{equation}
\label{39}
\Psi  = \int_0^{x_n} \D_{x_n} \Psi(t) dt =   x_n  I_1 (\D_{x_n} \Psi)  . 
\end{equation}

\medbreak

{\bf a) }  Suppose that  for some $\alpha \in [ 0, a [$ there is 
$\mu \in ]0, 1[$ such that 
\begin{equation}
\label{recreg}
 \Psi  \in  x_n^\alpha C^{ \mu }(\cD'_{\tg} ), \quad  
 x_n\D_{x_n}  \Psi  \in  x_n^\alpha C^{ \mu }(\cD'_{\tg} )
\end{equation}
  From \eqref{32} \eqref{33}, this is true when
$\alpha = 0 $ and $\mu  = 1/2$. 

  This condition implies that $ F  =     x_n^\alpha F_1 $ with $F_1 \in C^{ \mu }(\cD'_{\tg} )$.
  Therefore, $J_{a} (F ) = x^\alpha J_{ a -  \alpha }  (  F_1  ) $ and \eqref{37} implies that
\begin{equation*}
\D_{x_n} \Psi  =  x_n^a \theta  +   x_n^\alpha J_{a-\alpha} (F_1)  \in x^a C^{1/2}(\cD'_\tg)  + x_n^\alpha  C^{\mu} (\cD'_{\tg})
\subset x_n^\alpha C^{\mu_1} (\cD')
\end{equation*} 
for some $\mu_1 > 0$. 
  This shows that \eqref{recreg} holds with $\alpha$ replaced by $\alpha+1$, decreasing
$\mu$ if necessary.  

  Moreover, if \eqref{recreg} holds for some $\alpha$, it also holds for all $\alpha' < \alpha$,
with possibly smaller  $\mu$'s.
   Therefore, after a finite number of iterations, one obtains that \eqref{recreg}  is satisfied for
some $\alpha \in ]a, a+1 [$ and $\mu > 0$.

 \medbreak 

{\bf b ) } The property \eqref{recreg}  implies again that $F =  x_n^\alpha F_1 $ with $F_1 \in C^{\mu } (\cD'_{\tg})$.
 We now write
\begin{equation*}
\begin{aligned}
J_{a} ( F ) =  &   -  x_n^\alpha  I_{\alpha  -a}  (  F_1  )
                   + x_n ^a \int_0^\delta t^{ \alpha - a -1   } F_1(t) dt  \\
               & \in  x_n^\alpha  C^{\mu}(\cD'_y)  +  x_n^a \cD'_{\tg}  \subset x_n^a C^{\mu_1}(\cD'_y)
\end{aligned} 
\end{equation*}
for some $\mu_1 > 0$. Here we have used that $ \alpha > a$. 

  With \eqref{37} this implies that $\D_{x_n} \Psi = x_n^a  \Psi_1 $ with $\Psi_1  \in C^{\mu_1}(\cD'_\tg)$. Plugging this equation in \eqref{39}, implies that
\begin{equation}
\Psi = x_n^{a+1}  I_{a +1} (\Psi_1)  \in x_n^{a+1} C^{\mu_1}(\cD'_\tg).
\end{equation}

 Summing up, we have proved that $ \Phi = x_n^{- (a+1)} \Psi \in C^{\mu_1}(\cD'_\tg)$ and
$x_n \D_{x_n} \Phi = x_n^{-a} \D_{x_n} \Psi - (a+1) \Phi  \in   C^{\mu_1}(\cD'_\tg)$.

  \medbreak
{\bf c) } Inspecting the proof shows that the index $\mu $  and the Banach space
$E$ such that \eqref{36} holds are independent of $\Psi$ and $f$.
  Moreover, each step in the preceding proof can be converted into an estimate
ending with \eqref{estnorm}.
\end{proof}
 

\section{H\"older  estimates of $\Phi$. }
  We continue the proof of Theorem~\ref{mainestr}, proving the H\"older smoothness
of $\Phi$ and the estimate \eqref{mest1loc} of Theorem~\ref{mainestr}.
\begin{prop}
\label{regcmu}
  Suppose that  $\Psi$ and $f \in L^p ([0, \delta ] \times \omega) $, $p > n$,  are as in Theorem~$\ref{mainestr}$.
  Then for all relatively compact subset $\omega_1 \subset \omega$  and all $\mu \le  1-n/p$,
$\Phi = x_n^{-(a+1)} \Psi  \in C^\mu(\overline \omega_1 \times [0, \delta ] )$ and
$x_n\Phi \in C^{\mu+1}(\overline \omega_1 \times [0, \delta ] )$, with
\begin{equation}
\label{estregcmu}
 \Vert \Phi  \Vert_{ C^\mu(\overline \omega_1 \times [0, \delta ] )}  +
   \Vert x_n    \Phi \Vert_{C^{\mu+1}(\overline \omega_1 \times [0, \delta ] )}  
 \le   C \Big( \Vert f  \Vert_{L^p}
 + 
  \Vert \Psi \Vert_{H^1}    \Big)  . 
\end{equation} 
\end{prop}
  
  The equation  \eqref{eqloc}  for $\Psi$ implies that  $\Phi$ satisfies  
\begin{equation}
\label{geneq}
\cL \Phi :=  - x_n P_2 (x, \D_x) \Phi   -  Q_1(x, \D_x) \Phi   + R (x) \Phi  =  f , 
\end{equation}
   where 
  \begin{equation*}
  \begin{aligned}
  & P_2(x, \xi) = \sum_{j, k} p_{j, k} (x) \xi_j \xi_k, \\
 & Q_1(x, \xi)  =   2(a +1) \sum_j  p_{j, n} (x) \xi_j       - a  \sum_j p_j(x)  \xi_j    , \\  
& R(x) =  a (a+1) \big( p_n(x) -  p_{n, n} (x) \big) / x_n . 
\end{aligned}
\end{equation*}
 From \eqref{Zloc}, it follows that $R$ is smooth up to the boundary $ \{ x_n = 0 \}$
 and that 
 \begin{equation}
\label{Zloc2}
Q_1(x', 0, \xi) = ( a+2 ) \sum_{j=1}^n p_{j, n}(x', 0) \xi_j . 
\end{equation}
   Without loss of generality, we can also assume that $p_{n, n} = 1$. 
   
\bigbreak 
  Before proving the proposition above, we introduce  several  notations.
  For $\mu > 0$,  $ \mu \notin \NN  $, we denote by  $C^\mu (\overline \RR^n_+ ) $ the space of
bounded functions on $\RR^n_+$ which are uniformly H\"older continuous with exponent $\mu$.
  Next, we denote by $C^\mu_k(\overline \RR^n_+)$, $ k \in \NN$,  the space of functions
$u \in C^\mu (\overline \RR^n_+ ) $ such that $x_n^j u \in C^{\mu + j}$ for all $j \le k$.
$C^\mu_{\infty} $ is the intersection of the $C^\mu_k$.
 We refer to  \cite{BCM} for a definition of these spaces involving  a Littlewood-Paley 
analysis and their extension to negative $\mu$ and real $k$.
  We also refer to the same paper for a precise definition of the spaces $C^{\mu, \nu}_k(\overline \RR^n_+) $,
where the second index $\nu \in \RR $ measures the additional tangential smoothness.
  Finally, we introduce the spaces
   $C^{\mu, - \infty }_k = \bigcup_\nu  C^{\mu, \nu}_k$.

  Proposition~\ref{regcmu}  is a direct consequences of the results proved in \cite{BCM}, applied to second order equations \eqref{geneq}.
  There are three assumptions in \cite{BCM} to satisfy.

\begin{itemize}
\item  (H1)  The second order part $P_2$  is strongly elliptic.
\item  (H2) The indicial polynomial,
 $\lambda  P_2 (x', 0, 0 , 1)  +  Q_1(x', 0, 0, 1)$ has no roots in the strip
$\mu_1 \le \re \lambda \le \mu_2$.
\item (H3) The differential operator
  $L^0 = x_n P_2 (x', 0 ,  i \eta , \D_{x_n} )  + Q_1 (x', 0, i \eta,  \D_{x_n}) $ is 
injective in the space  $C^{\mu +1}_\infty (\overline \RR_+ )$ for all $\eta \in \RR^{n-1} \backslash \{0\}$
and all $x'$.
\end{itemize}
 
\bigskip

\noindent 1) For equation \eqref{geneq} the ellipticity assumption (H1) follows from \eqref{ellip}.

\smallbreak

\noindent 2)  By \eqref{Zloc2}, the indicial polynomial is $(\lambda + a+2)$.
    An important feature is that its roots are independent of $x'$.
    Since $a > 0$,  (H2)  is satisfied for all $-1 < \mu_1 \le \mu_2  $.
 
\noindent 3)  Using again \eqref{Zloc2},  we have
\begin{equation}
  L^0 = x_n  ( Z^2 -  \rho^2)   + (a+2) Z = Z (x_n Z) + (a+1) Z - x_n \rho^2, 
\end{equation}
 with 
\begin{equation}
\begin{aligned}
& Z =  \D_{x_n}  + i \sum_{j < n} \eta_j p_{j, n}(0, x'),
 \\
& \rho^2 = \sum_{j, k < n} p_{j, k}(x', 0) \eta_j \eta_k. 
\end{aligned}
\end{equation}

 Since $\rho >0$ when $\eta \ne 0$, the bounded functions $u$ in the kernel of $L^0$ are smooth and rapidly decreasing at infinity.
 Moreover, they are  smooth up to $0$, by the classical analysis of Fuchsian singularities.
 Thus, the following integration by parts are jusitified:
 \begin{equation*}
0 = \Re \int_0^\infty    L^0 u  \, \overline u  dx = - \int_{0}^\infty x_n
\big(  |  Z  u|^2 + \rho^2 u^2 \big) dx_n 
- \frac{a+1}{2} | u(0 |^2. 
\end{equation*}
 This implies that $u= 0$, hence that the spectral condition (H3) is satisfied if $\mu + 1 > 0$.

 Therefore we are in position to apply the results of \cite{BCM}.

\begin{proof}[Proof of Proposition~$\ref{regcmu}$]
  First, we note that,  because $p > n$, the Sobolev embedding theorem implies that the right-hand side of \eqref{geneq}  satisfies
\begin{equation}
\label{evid}
f  \in L^p \subset  C^{-n/p} (\omega \times [0, \delta] ).
\end{equation}

\noindent {\bf a) } $\Phi$ vanishes for $x_n \ge \delta/2$, so we can extend it by
 $0$ for  $x_n \ge \delta$. 
 By proposition~\ref{prop32}  we know that there are $\mu' \in ]0, 1/2]$ and a Banach space $E \subset \cD'(\omega)$ such that  $\Phi \in C^{\mu'}_1 (\overline \RR_+; E)$. 
  Therefore,  Proposition 2.3.2 of \cite{BCM}, implies that
 \begin{equation*}
\chi  \Phi \in C^{\mu' , - \infty}_1  (\overline \RR^n_+)
\end{equation*}
 if $\chi \in C^\infty (\omega)$ is equal to $1$ on $\omega_1$.
 This implies that there is a (large negative) integer $\nu $ such that
  \begin{equation}
 \label{estmu1}
\chi  \Phi \in C^{\mu_1 + 1, \nu }_1  (\overline \RR^n_+),  
\end{equation}
  if $\mu_1 = \mu' - 1  > -1$.

\medbreak
\noindent  {\bf b) }  Let $\chi_1 \in C^\infty_0(\omega)$, equal to  1 on $\omega_1$ and such
that $\chi = 1$ on a neighborhood of the support of $\chi_1$.
  There holds
  \begin{equation}
\label{eqloc21}
\cL \chi_1 \Phi = \chi_1 f  + [ \cL, \chi_1 ] \chi \Phi 
\end{equation}
  By \eqref{evid},
$\chi_1 f \in   C^{\mu - 1} (\overline \RR^n_+ ) \subset  C^{\mu - 1, 0} (\overline \RR^n_+ )$
with $\mu = 1 - n/p$.
  By \eqref{estmu1}, $[ \cL, \chi_1] \chi \Phi \in  C^{\mu_1 + 1, \nu } (\overline \RR^n_+ )$.
  Thus the right hand side of \eqref{eqloc21}  belongs to
$C^{ \mu_2, \nu} $  with $\mu_2 = \mu -1  < 0 < \mu_1 + 1 = \mu'$. 
  Since the assumption (H2) is satisfied for the pair $(\mu_1, \mu_2)$, Theorem 4.1 in \cite{BCM}
implies that
$\chi_1 \Phi \in C^{\mu_2 +1, - \infty}_1 $, that is
\begin{equation}
 \label{estmu2}
\chi_1  \Phi \in C^{\mu_2 + 1, \nu_1 }_1  (\overline \RR^n_+),  
\end{equation}
for some $\nu_1$.
    
\noindent {\bf c) }  Let $\chi_2 \in C^\infty_0(\omega)$, equal to  1 on
$\omega_1$ and such that   $\chi = 1$ on a neighborhood of the support of $ \chi_2$.
  There holds
\begin{equation}
\label{eqloc22}
\cL \chi_2 \Phi = \chi_2 f  + [ \cL, \chi_2 ] \chi_1 \Phi . 
\end{equation}
  By \eqref{evid},
$\chi_1 f \in   C^{\mu - 1} (\overline \RR^n_+ ) \subset  C^{\mu_2, 0} (\overline \RR^n_+ )$.
  By \eqref{estmu1}, $[ \cL, \chi_1] \chi \Phi \in  C^{\mu_2 + 1, \nu_1 } (\overline \RR^n_+ )
\subset C^{\mu_2 , \nu_1 +1  } (\overline \RR^n_+ )$. 
   By Theorem 5.2 in \cite{BCM}, we obtain that $\chi_2 \Phi \in C^{\mu_2 +1 , \nu_2}_1   $, with
$\nu_2 = \min (0, \nu_1 + 1)$.
   
   Therefore, after   finitely many iterations, we conclude that there is 
$\chi_* \in C^\infty_0(\omega)$, equal to  1 on $\omega_1$ such that
\begin{equation}
 \label{estmu3}
\chi_*  \Phi \in C^{\mu_2 + 1, 0 }_1  (\overline \RR^n_+).   
\end{equation}

\noindent {\bf d)} Finally, introducing $\chi_\flat  \in C^\infty_0(\omega)$, equal to 1
on $\omega_1$ and such that $\chi_*= 1$ on a neighborhood of the support of
$\chi_\flat $ and writing the equation for $\chi_\flat \Phi$,
Theorem 4.2 in \cite{BCM} implies that
$\chi_\flat \Phi \in C^{\mu_2 + 1 }_1 (\overline \RR^n_+) = C^{\mu }_1 (\overline \RR^n_+)$
implying the proposition.
\end{proof}


\section{$L^p$ estimates}
  In this section, we finish the proof of Theorem~\ref{mainestr}, proving the
$L^p$  estimates \eqref{mest2loc}.
  We recall that regularity result is performed in the general setting $\Omega\subset \RR^n$ with $n\ge 1$.
 Given $\omega_1 $ as in the theorem, consider $\omega_2$ relatively compact
in $\omega$ such that $\overline \omega_1 \subset \omega_2$. 
 By Proposition~\ref{regcmu}, $\Phi$ is of class $C^\mu $ on $\omega_2 \times [0, \delta]$ for $\mu \le  1- n/p $. Consider $\chi \in C^\infty_0(\omega_2)$ with $\chi = 1$ on $\omega_1$.
 Let $\phi  = \chi \Phi \in C^\mu_1  (\overline \RR^n_+)$. Then
\begin{equation}
\label{eqterm}
\cL \phi = g  :=  \chi f + [\cL, \chi] \Phi \in L^p (  \RR^n_+).
\end{equation}
   Moreover, by Proposition~\ref{regcmu}, choosing $p_0 > n$ and $\mu_0 = 1 - n/p_0$, 
   there holds for $f \in L^p $,  $p \ge p_0$: 
\begin{eqnarray}
\label{reginterm}
\Vert \phi \Vert_{C^{\mu_0}_1(\overline \RR^n_+)}  
& \le & C \big( \Vert f \Vert_{L^{p_0}}  + \Vert \Psi \Vert_{H^1} \big),
\\  
\Vert g  \Vert_{L^p (  \RR^n_+)}  
&\le & C \big( \Vert f \Vert_{L^{p}}  + \Vert \Phi \Vert_{C^{\mu_0}_1(\overline \RR^n_+)} \big). 
\end{eqnarray}
 Therefore, the $L^p$ estimates \eqref{mest2loc} follow from the next result.
 \begin{theo}
\label{thmainestLp} Let $\mu \in ]0, 1[$. 
Suppose that $\phi \in C_1^\mu (\overline \RR^n_+)$  
 has compact support in $\omega \times [0, \delta[ $  and 
 $\cL  \phi  = g  \in L^p(\RR^n_+)$.
Then 
 $\phi  \in   W^{1, p}_{1, loc}(\overline \RR^n_+)$.

Moreover, there is a constant $C$ such that for all such $\phi$ and $p \in [2, \infty [$
  \begin{equation} 
\label{mainestx}
\Vert  x_n    \D_j \D_k  \phi   \Vert_{L^p(\RR^n_+)}  + 
\Vert  \D_j   \phi   \Vert_{L^p(\RR^n_+)}  \le    C  \Big(  p     \Vert   g  \Vert_{L^p(\RR^n_+) }
+ \Vert \phi  \Vert_{C^\mu_1} \Big) . 
\end{equation}
\end{theo}
  
\subsection{Preliminary results}
 In this section we recall some known results about Calderon-Zygmund operators (\cite{St}).
We consider operators $T$ acting in $L^\infty_{comp} (\overline \RR^n_+)$ 
 (the space of bounded functions with compact support) with
kernel $T(x,y)$ locally integrable away from the diagonal $\{ x = y \}$.

\begin{prop}
 \label{action1} Suppose that the kernel $K(x,y)$ satisfies on $\RR^n_+ \times \RR^n_+$: 
\begin{equation}
\label{maj1}
| K(x, y ) | \le \frac{C}{| x - y|^{n-1}}, \quad  | \D_x K(x, y ) | \le \frac{C}{| x - y|^{n}}.
\end{equation}
Then, the operator 
\begin{equation*}
T f(x) = \int K(x, y) f(y) dy
\end{equation*}
acts from $L_{comp}^p(\overline \RR^n_+)$  to $C^\mu(\overline \RR^n_+)$ for all $\mu < 1-n/p$.
\end{prop}
 
 \begin{proof}
  The first estimate for $K$ implies that $T$ maps $L_{comp}^\infty(\overline \RR^n_+)$ to
$L^\infty$.
  Moreover, the second estimate implies that
\begin{equation*}
| K(x, y) - K(x', y) | \le C \frac{ | x - x' |} { | x - y |^{n}}  \quad \mathrm{for}
\ \ | x - x' | \le \frac{1}{2} | x - y |. 
\end{equation*}
Thus, interpolating with the first estimate yields 
\begin{equation}
\label{maj1b}
| K(x, y) - K(x', y) | \le C \frac{ | x - x' |^\mu } { | x - y |^{n-1 + \mu}}  \quad \mathrm{for}
\ \ | x - x' | \le \frac{1}{2} | x - y |. 
\end{equation}
 We write
\begin{equation*}
\begin{aligned}
T f(x) - Tf(x') =&  \int_{| x - y | \ge 2 | x - x' |}  \big( K(x, y) - K(x', y)  \big) f(y) dy 
\\  & + \int_{| x - y | \le  2 | x - x' |}  \big( K(x, y) - K(x', y)  \big) f(y) dy . 
\end{aligned}
\end{equation*}
By \eqref{maj1b} the first integral is $O(| x - x' |^\mu)$. Note that for 
$| x - y | \le  2 | x - x' |$, there also holds 
$| x' - y | \le  3 | x - x' |$. Therefore, the second integral is 
\begin{equation*}
O \Big( \int_{ | y | \le 3 | x - x'|}  \frac{d y} {| y|^{n-1}} \Big)  = O\big( | x - x' | \big) . 
\end{equation*}
 \end{proof}

 \begin{prop}
 \label{action2}
  Suppose that $T$ is a bounded operator in 
  $L^2(\RR^n_+)$ with  kernel $K(x,y)$ satisfying for $x \ne y$ : 
 \begin{equation}
\label{maj2}
| K(x, y ) | \le \frac{C}{| x - y|^{n}}, \quad  | \D_x K(x, y ) | \le \frac{C}{| x - y|^{n+1}}.
\end{equation}
   Then  $T$ maps $L^p(\RR^n_+)$ to $L^p(\RR^n_+)$ with norm $O(p)$ for
all $p \in [2, + \infty[$. 
 \end{prop}
 
\begin{proof}
 The adjoint operator $T^*$ is bounded in $L^2$ and its kernel
 $K^*(x, y) = \overline K(y, x)$ which therefore satisfies
\begin{equation}
\label{maj3}
| K^*(x, y ) | \le \frac{C}{| x - y|^{n}}, \quad  | \D_y K^*(x, y ) | \le \frac{C}{| x - y|^{n+1}}.
\end{equation}

  The Calderon-Zygmund theory implies that $T^*$ is bounded from $L^1$ to the space
weak-$L^1$ (see e.g. \cite{St}).
  Therefore, by Marcinkievic's interpolation theorem, $T^*$ maps $L^p$ to $L^p$ for $p \in ]1, 2]$
with norm $O (1/(p-1)$ (same reference).
  By duality, this implies that $T$ maps $L^p$ to $L^p$ for $p \in [ 2, \infty [ $ with norm
$O (p)$.
\end{proof}

\subsection{Parametrices}
 Recall that
 $\cL \approx -  x_n P_2  - Q_1$, see \eqref{geneq}.
 For $y \in \overline \omega$, we denote by $\cL_y$ the operator 
\begin{equation}
\label{56f}
\cL_y (x, \D_x) :=  -  x_n  P_2 (y, \D_x )       -      (a+2)  \big(   \D_{x_n}    
+ \sum_{j < n}  p_{j, n}(y)  \D_j \big) 
\end{equation}
As above, we have assumed as we may that $p_{n, n} = 1$. 
For a given $y$, there is a linear transformation 
\begin{equation}
\label{58}
\tilde x =  T(y) x , \quad  \mathrm{with} \ \    \tilde x_n = x_n 
\end{equation}
such that,  in these variables  $\cL_y$  is  transformed to 
\begin{equation}
\label{56m}
\widetilde  \cL  =  -  \tilde x_n \Delta_{\tilde x}  - ( a + 2)  \D_{\tilde x_n}. 
\end{equation}
 
 According to \cite{GS}, Lemme 1,  the  fundamental solution  of 
 $ \widetilde  \cL $ is 
 \begin{equation}
\label{gf2}
\widetilde E  (\tilde x, \tilde y   )  = 
\int_0^1 \tilde F  (\tilde x,  \tilde y  , \theta) d\theta,
\end{equation}
with 
\begin{equation} 
\begin{aligned} 
 \tilde F  (\tilde x,  \tilde y  , \theta)  &  =  \gamma \ (\tilde y_n)^{ a+1}  
  A ^{ - (a  + n )} 
 \left( \theta   (1- \theta)  \right)^{ a  /2}  ,\\
 A^2 (\tilde x, \tilde y, \theta) & =   \theta D^2 + (1- \theta) \check D^2. 
\end{aligned}
\end{equation}
with $\gamma $ some constant depending on $a$ and  
$$
D^2 = | \tilde x_n - \tilde y_n |^2 + | \tilde x' - \tilde y'  |^2   , \quad 
 \check D^2 = | \tilde x_n + \tilde y_n |^2 + | \tilde x' - \tilde y'  |^2  .
$$
More precisely, for $\eps \in ]0, 1[$, let 
 \begin{equation}
\label{gfeps}
\widetilde E^\eps   (\tilde x, \tilde y   )  = 
\int_0^{1- \eps}  \tilde F  (\tilde x,  \tilde y  , \theta) d\theta. 
\end{equation}
This function is singular only on $\tilde y_n = 0$ and there holds
\begin{equation}
\label{solf1}
\widetilde L (\tilde x, \D_{\tilde x}) \tilde E^\eps(\tilde x, \tilde y)  =  
 \widetilde G^\eps  (\tilde x, \tilde y)  , 
\end{equation}
where 
\begin{equation}
\label{errsolf}
 \tilde G^\eps  (\tilde x,  \tilde y )     =  2 (n - 2 + a ) \gamma \ \tilde y_n^{ a+2}  
  A ^{ - (a +2 + n )} 
 \big( \eps   (1- \eps)  \big)^{(a+2) /2} . 
\end{equation}
According to \cite{GS} , Appendix A, there holds :

\begin{lem}
\label{lemapproxid}
$  \tilde G^\eps   (\tilde x,  \tilde y )  $ is an  approximation of the 
identity as $\eps $ tends to zero: it is nonnegative, converges uniformly to 
$0$ on compacts of $ \overline \RR^d_+ \times \overline \RR^d_+ \setminus \{ \tilde x = \tilde y\} $ and 
for all bounded open set $\Omega \subset \RR^n_+$ and all $\tilde x \in \Omega$
\begin{equation*}
\lim_{\eps \to 0}  \int_{\Omega}    \tilde G^\eps   (\tilde x,  \tilde y )  d\tilde y
=  1 . 
\end{equation*}
\end{lem}
 From $\widetilde E$, we derive fundamental solutions for the operators $\cL_{z}$:
 \begin{equation}
\label{fundsol}
E_z (x, y) =  | \det T'(z) | \, E \big(T(z) x, T(z) y\big)
\end{equation}
and their approximate versions
\begin{equation}
\label{fundsolapp}
E^\eps_z (x, y) =  | \det T'(z) | \, E^\eps \big(T(z) x, T(z) y\big)
\end{equation}
Finally, we define the parametrices
 \begin{equation}
\label{param}
E (x, y)  = E_y (x, y) , \quad 
E^\eps (x, y)  = E^\eps_y (x, y) . 
\end{equation}
Similarly, we define the $G^\eps_z(x, y)$ and $G^\eps(x, y)$.

\begin{prop}
\label{propE} 
Suppose that  $a \ge 1$. For $k \in \NN$,  there  is a  constant  $C_k$ such that for  $(x, y) \in \RR^n_+ \times \omega $, $x \ne y$, $\eps \in ]0, 1[$, 
there holds
\begin{equation}
\label{esE}
  | \na_x^k E^\eps (x, y) | \le \frac{C_k}{ | x - y |^{ n+k -1}}, 
  \end{equation}
  and for $k \ge 1$
  \begin{equation}
  | x_n \na_x^k E^\eps (x, y) | \le \frac{C_k}{ | x - y |^{ n+k -2}}. 
\end{equation}
The same estimates hold for $E$.  
\end{prop} 

\begin{proof}
It is sufficient to prove the estimates for 
$E^\eps_z(x, y)$ and therefore for $\widetilde E(\tilde x, \tilde y)$. 
For simplicity, we drop the tildes in the proof below. 
On $\RR^n_+ \times \RR^n_+$, there holds 
\begin{equation}
\label{A}
A  \ge | x -  y |  , \quad   A  \ge (1 - \theta)^{1/2} (x_n+ y_n). 
\end{equation}
Hence, for $a \ge - 1$, 
$$
  y_n^{ a+ 1}  
  A ^{ - (a   +n )}  \le | x - y |^{1 - n} 
   (1- \theta)^{ - (a +1) /2}  , 
 $$
 implying 
\begin{equation*}
| F(x, y, \theta) |   \le \gamma | x - y |^{ 1 - n} \theta^{ a/2} (1- \theta)^{-1/2} . 
\end{equation*}
Thus, integrating in $\theta$, 
\begin{equation}
\label{EstE}
|E^\eps (x, y)| \le C | x - y |^{1-n}
 \end{equation} 

\noindent { \it Estimates on the derivatives of $E$.}
   Differentiating  $k$ times  $A^2$, one proves by induction that
\begin{equation}
| \na^k _x A | \le C_k  A^{1 - k } .
\end{equation} 
Thus 
\begin{equation}
| \na_x^k  F | \le  C_k  {y_n}^{ a+1 }   A ^{ - (a  + k + n)} 
 \left( \theta   (1- \theta)  \right)^{a /2} . 
 \end{equation} 
 Using \eqref{A}  yields 
$$
|  \na^k_x F(x,y)    | \le  C     | x - y |^{   1 - k - n}  
  \theta^{a/2}    (1- \theta) ^{ -  1 /2} ,  
$$
implying 
$$
|  \na^k_x E^\eps (x, y)      | \le  C     | x - y |^{ 1 - k  - n} .   
$$

\noindent {\it Sharper estimates of $x_n \na^k_x E $.}
 There holds
 \begin{equation}
\label{nA}
A^2 = | x - y |^2  + 4 (1- \theta) x_n y_n. 
\end{equation}

 Let us assume first that
 \begin{equation}
\label{case1}
 x_n y_n \le  2  | x - y |^2.
\end{equation}
 Then,
\begin{equation}
\label{case1b}
( x_n + y_n)^2 \le (x_n - y_n)^2 + 4 x_n y_n \le  C | x - y |^2. 
 \end{equation}
 In this case, we use that $A \ge | x - y|$, without loosing any information and 
 \begin{equation*}
| x_n \na_x^k F(x, y) | \le C x_n y_n^{a+1} | x - y |^{ - (n + k + a)} 
(\theta (1 - \theta))^{a/2} 
\end{equation*}
 Integrating in $\theta$ yields 
  \begin{equation*}
| x_n \na_x^k E(x, y) | \le C x_n y_n^{a+1} | x - y |^{ - (n + k + a)}  . 
\end{equation*}
With \eqref{case1b}, $x_n y_n^{a+1} \le C | x - y |^{a+2} $, therefore
\begin{equation}
\label{case1E}
| x_n \na_x^k E^\eps (x, y) | \le C   | x - y |^{ - (n + k  - 2)}.   
\end{equation}

Consider next the case
 \begin{equation}
\label{case2}
 x_n y_n \ge  2  | x - y |^2. 
\end{equation}

In this case,  
\begin{equation}
\label{case2b}
 5 x_n y_n \ge   2  x_n^2  + 2 y_n^2   \quad  \Rightarrow \quad
 \frac{1}{C} \le   \frac{x_n}{y_n} \le C . 
\end{equation}
Let 
\begin{equation*}
\rho := \frac{ | x - y|^2}{4 x_n y_n } \in ]0, \frac{1}{2}]. 
\end{equation*}
We use the estimates 
\begin{equation*}
 \begin{aligned}
 &  A^2 \ge |  x - y |^2 \quad  & \mathrm{when} \quad  1 - \theta  \le \rho, 
 \\
&  A^2 \ge 4 (1- \theta) x_n y_n  \quad  & \mathrm{when} \quad 1 -  \theta  \ge \rho.
\end{aligned}
\end{equation*}
Therefore, 
\begin{equation*}
\begin{aligned}
| x_n \na_x^k E^\eps (x, y) | \le C  &   \frac{ x_n y_n^{a+ 1} } { | x - y |^{ n + k + a }} 
\int_{1-\rho}^{1 }  \big( \theta (1 - \theta) \big)^{a/2} d \theta \\
& +    C  \frac{ x_n y_n^{a+1} } { (x_n y_n)^{ (n + k + a)/2}} 
\int_0^{1-\rho}  \theta^{a/2}  (1 - \theta)^{- (n+k)/2} d \theta . 
\end{aligned}
\end{equation*}
Thus 
\begin{equation*}
\begin{aligned}
| x_n \na_x^k E^\eps (x, y) | \le C  &   \frac{ x_n y_n^{a+1} } { | x - y |^{ n + k + a}} 
\rho^{a/2  }   \\
& +    C  \frac{ x_n y_n^{a+1} } { (x_n y_n)^{ (n + k + a)/2}} 
\rho^{1- (n+k)/2} .
\end{aligned}
\end{equation*}
For  the last estimate we have used that $n+k > 2$. Therefore, 
\begin{equation*}
\begin{aligned}
| x_n \na_x^k E^\eps(x, y) | \le C  &   \frac{ x_n y_n^{a+1} } { | x - y |^{ n + k + a}} 
\frac{ | x - y |^{ a +2 } }{ (x_n y_n)^{(a+2)/2  }}    \\
& +    C  \frac{ x_n y_n^{a+ 1} } { (x_n y_n)^{ (n + k + a)/2}} 
\frac{ (x_n y_n)^{(n+k -2)/2 } }{| x - y |^{ (n+k -2)}} 
\end{aligned}
\end{equation*}
\begin{equation}
| x_n \na_x^k E^\eps (x, y) | \le C   \    \frac{   y_n^{a/2} }  { x_n^{a/2  }}\  
\frac{1}{ | x - y |^{ n + k - 2 }} . 
   \end{equation}
Together with \eqref{case2b} this implies 
\begin{equation}
| x_n \na_x^k E^\eps (x, y) | \le C   \  
\frac{1}{ | x - y |^{ n + k - 2 }} . 
   \end{equation}
This finishes the proof of the proposition.
\end{proof}

\begin{lem}
\label{lemestK}
The following identity holds 
\begin{equation}
\label{comp1}
\cL  (x, \D_x) E^\eps (x, y) =   G^\eps(x, y)  +  K^\eps (x, y) 
\end{equation}
where the  kernels $K^\eps$  satisfy 
uniformly in $\eps \in ]0, 1[$: 
\begin{equation}
\label{Ke}
| K^\eps(x, y ) | \le \frac{C}{| x - y|^{n-1}}, \quad  | \D_x K^\eps(x, y ) | \le \frac{C}{| x - y|^{n}}.
\end{equation}
Moreover, the kernels $K^\eps$ converge uniformly on compacts of 
$\{ x \ne y\}$ 
to $K (x, y)$ which satisfy  $\eqref{Ke}$  outside the diagonal. 

\end{lem}

\begin{proof} 
The definition of $E^\eps$ implies that it is singular only on $y_n = 0$. 
Moreover, 
\begin{equation*}
\cL  (x, \D_x) E^\eps (x, y) =   
\cL_y (x, \D_x) E^\eps (x, y)  + K^\eps(x, y) = G^\eps(x, y)  +  K^\eps (x, y)
\end{equation*}
where 
\begin{equation}
\label{defK}
\begin{aligned}
K^\eps(x, y) = x_n \sum& \big( p_{j, k} (y) - p_{j,k}(x)\big) \D_{x_j} \D_{x_k}  E^\eps \\
 & + \sum  \big(  ( a+ 2)  p_{j, n}(y)  - q_j(x) \big) \D_{x_j} E^\eps   +  R(x) E^\eps . 
\end{aligned} 
\end{equation}
where the $q_j$ are the coefficients of $Q_1$.
  In the first sum, the coefficient of $ x_n \D_{x_j} \D_{x_k}  E^\eps$ is
$O (| x - y |)$.
 By \eqref{Zloc2}, the coefficient of $\D_{x_j} E^\eps$ in the second sum is
$O( x_n ) +O( | x - y |)$.
  Together  with Proposition \ref{propE}, this implies the estimates \eqref{Ke}.

 From \eqref{defK}, it is clear that $K^\eps$ converges to a kernel $K$ which satisfies 
 \eqref{Ke}. 
\end{proof} 

\begin{lem}
\label{lemapproxid2}
The kernels $G^\eps(x, y)$ are non negative, converge to $0$ as $\eps$ tends to $0$ uniformly
on compacts subsets of
$(\overline \RR^n_+ \times \overline \omega ) \setminus \{ x = y \}$.
  Moreover the integrals $$\int_{\omega \times ]0, \delta[ }  G^\eps (x, y) d y $$
are uniformly bounded for $\eps \in ]0, 1[$ and $x$ in a compact subset of
$\omega \times ]0, \delta[  $, and converge to $1$ as $\eps$ tends to $0$.
\end{lem}
\begin{proof}
By \eqref{errsolf}, \eqref{nA} and  \eqref{58}, it follows that 
\begin{equation*}
\begin{aligned}
G^\eps_z(x, y) &  = | \det T'(z) | \widetilde G^\eps \big(T(z) x, T(z) y) \big) \\
 &  =  | \det T'(z) |  \frac{  y_n^{a+2}  
\big( \eps (1- \eps) \big)^{(a+2)/2} }{( | T(z) (x - y) |^2 + 4 \eps x_n y_n)^{(a+n+2)/2} }
\end{aligned}
\end{equation*}
is nonnegative.  
Thus
\begin{equation}
\label{estG}
 G^\eps_z(x, y)  + | \D_z G^\eps_z(x, y) | \le  \frac{ C y_n^{a+2}  \eps^{(a+2)/2}} {( | x - y |^2 + \eps x_n y_n)^{ (a + n+2)2}} . 
\end{equation}
In particular, $G^\eps(x, y) = G^\eps_y(x, y)$ is also dominated 
by the same bound. 
implying     that $G^\eps(x, y) \to 0$ as $\eps \to 0$  when 
$x \ne y$, uniformly on compact subsets.

 Integrating in the tangential variables first implies that
\begin{equation}
\label{estGxx}
\int_{\omega \times ]0, \delta[ } G^\eps(x, y) dy \le 
\int_0^L   \frac{ C y_n^{a+2}   \eps^{(a+2)/2}} {( | x_n - y_n |^2 + \eps x_n y_n)^{ (a + 3)/2}} 
dy_n. 
\end{equation}
 The integral  over  $\{ | x_n - y_n | \ge x_n/2$ tends to zero, uniformly 
 for $x_n $ in a compact subset of $] 0, \delta ] $. 
 The integral over the remaining interval  $\{ | x_n - y_n | \le x_n/2$ is 
 bounded by 
 \begin{equation*}
 \begin{aligned}
 \int_{\frac{x_n}{2}}^{\frac{3 x_n}{2} } 
    \frac{ C x_n^{a+2}  \eps^{(a+2)/2}} {( | x_n - y_n |^2 + \eps x_n^2)^{ (a + 3)2}}  dy_n 
   & =  \int_{- \frac{1}{2}}^{\frac{1}{2} } 
    \frac{ C   \eps^{(a+2)/2}} {( | y_n |^2 + \eps )^{ (a + 3)2}}  dy_n \\
   &  \le 
     \int_{\RR } 
    \frac{ C    } {( | y_n |^2 + 1  )^{ (a + 3)2}}  dy_n  < + \infty. 
\end{aligned}
\end{equation*}

Moreover, by \eqref{estG}  there holds
\begin{equation*}
 |   G^\eps_x(x, y) - G^\eps_y(x, y)  | \le  \frac{ C | x - y| y_n^{a+2}  \eps^{(a+2)/2}} {( | x - y |^2 + \eps x_n y_n)^{ (a + n+2)2}} . 
 \end{equation*}
 This implies that   
 \begin{equation*}
\lim_{\eps \to 0} \, \int_{\omega \times ]0, \delta[ } |   G^\eps_x(x, y) - G^\eps_y(x, y)  |   d y    = 0. 
\end{equation*}
 By  Lemma~\ref{lemapproxid} (see \cite{GS}), for $x \in \Omega$ there holds
  \begin{equation*}
\lim_{\eps \to 0} \, \int_{\omega \times ]0, \delta[ }    G^\eps_x(x, y)    d y    =  1,  
\end{equation*}
implying that 
  \begin{equation*}
\lim_{\eps \to 0} \, \int_{\omega \times ]0, \delta[ }    G^\eps (x, y)    d y    =  1.  
\end{equation*}
This finishes the proof of the lemma. 
\end{proof}

  We denote by $\bE$, $\bE^\eps$ etc  the operator with kernel
$E$, $E^\eps$ etc.
 In particular, since the kernels are smooth
for $y_n > 0$ and bounded on $\RR^n_+ \times \RR^n_+$, the operators
$\bE^\eps$, $\bG^\eps$ and $\bK^\eps$ are defined from  the space of integrable
functions with  support in $\omega \times [0, \delta]$  to the space of $C^\infty(\overline \RR^n_+)$.
 The lemma above implies that
\begin{equation}
\label{comp2}
\cL   \bE^\eps =  \bG^\eps + \bK^\eps . 
\end{equation}
  Moreover, the estimates \eqref{esE} and \eqref{Ke} imply that the operators
$\bE^\eps$ [resp.  $\bK^\eps$]  are uniformly bounded from $L^p (\omega \times ]0, \delta[)$
to $L^p_{loc}(\overline \RR^n_+)$ and converge in the strong topology to the operators
$\bE$ [resp. $\bK$]  defined by the kernels $E$ [resp. $K$].
    
\begin{prop}
\label{propactionE}
  The operator $\bE$ maps $L^p (\omega \times ]0, \delta[ )$ to
$  W^{1, p}_{1, loc}(\overline \RR^n_+)$.
  For all relatively compact open set $\omega_1 \subset \omega$, all
$\delta' < \delta$  and for all $\chi \in C^\infty_0 (\omega \times [0, \delta[ )$ such that
$\chi = 1$ on a neighborhood of $\overline \omega_1 \times [0, \delta]$ , there is
a constant $C$ such that for all $p \in [2, + \infty[$ and all $ g  \in L^p$ supported in
$\omega_1 \times ]0, \delta'[$, there holds
\begin{equation}
\label{mainestn}
\Vert  x_n \chi  \D_j \D_k  \bE g  \Vert_{L^p(\RR^n_+)}  + 
\Vert \chi    \D_j    \bE g \Vert_{L^p(\RR^n_+)}  \le   p  C \Vert   \bE g \Vert_{L^p(\RR^n_+)} .
\end{equation}

Moreover, for all  $p > n$ and $g \in L^{p}$  supported in
  $\omega_1 \times [0, \delta']$ and all  $\mu <  1 - n/p $, there holds  
\begin{equation}
\label{mainestn2}
\cL \bE  g  -    g  =  \bK g \in C^{\mu}_{1, loc} (\overline \RR^n_+) . 
\end{equation}
  \end{prop} 

\begin{proof}
  By Propositions~\ref{propE}, for $g \in L^p(\omega \times ]0, \delta[) $,
$\bE^\eps g $ and $\bK^\eps g $  are uniformly bounded in $L^p (\overline \RR^n_+) $ and converge
in $L^p_{loc}$ to $\bE g$ and $\bK g $ respectively.

  Suppose that $g$ is continuous with compact support in $\omega \times ]0, \delta [$. Then, by Lemma~\ref{lemapproxid2}, $\bG^\eps g $ converges to $g$ as $\eps$ tends to $0$.
This shows that the identity \eqref{mainestn2} is satisfied, in the sense of distributions, for
$g \in C^0 $ with compact support in $\omega \times ]0, \delta[$.
  By density, this identity extends  to $g \in L^2 (\omega \times ]0, \delta[)$.
  In particular, $\bE g \in  L^2_{loc}(\overline \RR^n_+)$ and
$\cL \bE g \in  L^2_{loc}(\overline \RR^n_+)$.
  Using the regularity properties of $L$ proved in \cite{BC}, this implies that $\bE$ maps  $L^2 (\omega \times ]0, \delta[)$ into the space  $H^{1}_{1, loc}(\omega \times [0, \delta[)$  of functions $u \in H^1_{loc}(\omega \times [0, \delta[)$ such that $x_n u \in H^2_{loc}(\omega \times [0, \delta[)$.

 In particular, this proves that for $\chi_1$ and $\chi_2 $ in $C^\infty_0(\overline \RR^n_+)$,
the operators $\chi_1 x_n \D_j \D_k \bE \chi_2 $ and $\chi_1  \D_j   \bE \chi_2 $ are bounded from
$L^2(\RR^n_+)$ to $L^2(\RR^n_+)$.
 Together with the estimates of Proposition~\ref{propE}
and with Proposition~\ref{action2}, this implies that $\bE$ satisfies the properties listed in the proposition and the estimates \eqref{mainestn}.

 The identity \eqref{mainestn2} holds for $g \in L^p$ thus for
$g\in L^\infty$.
 The properties of $\bK $ follow from Proposition~\ref{action1} together with the
estimates of Lemma~\ref{lemestK}.
\end{proof}

\subsection{Proof of Theorem~$\ref{thmainestLp}$}
 By assumption, $\phi $ and $g = \cL \phi$ have  compact support in
$\omega \times [0, \delta[$.
 Let $\chi $ and $\tilde \chi$  in $C^\infty_0( \omega \times [0, \delta [)$ 
such that $\chi = 1$ on a neighborhood of the support of $\phi$ and
$\tilde \chi = 1$ on a neighborhood of the support of $\chi$.

 By Propositions~\ref{propE} and \ref{action1}, $\tilde \chi \bE g \in  C^\mu_1(\overline \RR^n_+) $.
Let $\psi  = \chi \bE g $, which satisfies the
estimates \eqref{mainestn}.
 Moreover,
$$
\cL  \psi   = \chi \cL \bE g    + [\cL, \chi]    \bE   g   =  g   + \chi \bK  g   + 
[\cL, \chi]  \tilde \chi \bE g  .
$$ 
 In particular,
\begin{equation}
\label{equerr}
\cL (\phi  - \psi )  = \chi \bK g  +  [\cL, \chi] \tilde \chi \bE g .
\end{equation}
 By Proposition~\ref{propactionE}, $\chi \bK g \in C^\mu_1( \overline \RR^n_+)$ and
since $\tilde \chi \bE g \in  C^\mu_1(\overline \RR^n_+) $, the commutator 
$[\cL, \chi] \tilde \chi \bE g$ also belongs to  $C^\mu( \overline \RR^n_+)$. 
Thus the right hand side $h$  of \eqref{equerr} is in $C^\mu( \overline \RR^n_+)$

  We are now in position to apply the results of \cite{GS} (or \cite{BCM}), which imply that
\begin{equation*}
 \phi  - \psi  \in C^{\mu + 1}_{loc}(\overline \RR^n_+) ,  \quad
x_n(\phi - \psi) \in C^{\mu + 2}_{loc}(\overline \RR^n_+).
\end{equation*}
  Therefore,
\begin{equation}
\label{errxxx}
\begin{aligned}
\Vert  x_n    \D_j \D_k   (\phi - \psi)    \Vert_{L^\infty(\RR^n_+)}  + 
 & \Vert      \D_j (  \phi - \psi)   \Vert_{L^\infty (\RR^n_+)}  \le  \\
&    C  \Big(      \Vert   h  \Vert_{C^\mu(\overline \RR^n_+) } 
+ \Vert  (\phi - \psi)   \Vert_{C^\mu_1(\overline \RR^n_+)} \Big)
\end{aligned}
\end{equation}
with $h$ the right-hand side of \eqref{equerr}.

Together with the estimates of Proposition~\ref{propactionE} for the 
$L^p$ norm of 
$  x_n    \D_j \D_k  \psi  $ and 
$   \D_j   \psi$, this implies the theorem.


\section{Proof of global existence, uniqueness and regularity theorem.}

\noindent{\bf Regularity of weak solutions.}
  Suppose that $(v, \omega) $ is a weak solution in the sense of Definition~\ref{def21}.  
  Then
$v\in L^\infty(0,T;W^{1,p}(\Omega))$ for all $p<+\infty$ from the
elliptic regularity theorem.
 Therefore,  using
Di-Perna and Lions renormalization technics for transport equation \eqref{transp},
implies that   $\omega\in C(0,T; L^r(\Omega))$,  hence  that
$v\in C(0,T; W^{1,r}(\Omega))$ for all $r<+\infty$.
 The other properties of $v$   also follow from Theorem \ref{mainest}.

\bigbreak

\noindent{\bf Existence of a global weak solution.}
  We construct a global weak solution as the inviscid limit of
solutions to a Navier-Stokes system with artificial viscosity and boundary conditions.

\noindent{\it Global existence for a viscous vorticity-stream function formulation.}
  We define $b_\varepsilon = b + \varepsilon$ and consider the following system
\begin{equation}
\label{VSviscous}
\left\{
\begin{aligned}
& \partial_t (b^\varepsilon \omega^\varepsilon) + b v^\varepsilon\cdot\nabla \omega^\varepsilon
   - \varepsilon {\rm div}(b^\varepsilon\nabla \omega^\varepsilon) = 0 \hbox{ in } \Omega, \\
& b^\varepsilon\omega^\varepsilon|_{t=0} = b^\varepsilon\omega_0 \hbox{ in } \Omega,
  \qquad \omega^\varepsilon|_{\partial \Omega} = 0, \\
& - {\rm div}(\frac{1}{b}\nabla \Psi^\varepsilon) = b^\varepsilon\omega^\varepsilon\hbox{ in } \Omega,
    \qquad \Psi^\varepsilon|_{\partial\Omega} = 0, \\
& v = \frac{1}{b} \nabla^\bot \Psi^\varepsilon.\\
\end{aligned}\right.
\end{equation}
  The existence of a global weak solution $(v^\varepsilon,\omega^\varepsilon)$ for
$\omega_0 \in L^2(\Omega)$ follows standard techniques for equations of Navier-Stokes type
since $b^\varepsilon$ is a strictly positive function, see \cite{LOT}.
  From the a-priori estimates, we get the following uniform bounds with respect to $\varepsilon$
$$ \sqrt b^\varepsilon\omega^\varepsilon \in L^\infty(0,T;L^2(\Omega)), \qquad
   \sqrt b^\varepsilon\nabla\omega^\varepsilon \in L^2 (0,T; L^2(\Omega)).$$
 Moreover
$\sqrt b v^\varepsilon \in L^\infty(0,T;L^2(\Omega))$ uniformly.
\medbreak

\noindent{\it The viscous limit $\varepsilon\to 0$.}
 Let us now assume that $\omega_0\in L^\infty(\Omega)$.
 Since $\omega^\varepsilon \in L^\infty(0,T;L^2(\Omega))\cap L^2(0,T; H^1_0(\Omega))$, using the elliptic regularity for degenerate elliptic equation, we get $\partial_t \omega^\varepsilon \in L^\infty(0,T; H^{-1}(\Omega))$ and thus $\omega^\varepsilon \in C([0,T); L^2(\Omega)).$
 We multiply the viscous equation by
$|T_R(\omega^\varepsilon)|^{p-2} T_R(\omega^\varepsilon)$ where
$T_R(\omega^\varepsilon) = {\rm max}({\rm min}(\omega^\varepsilon,R),-R)$
for $R>0$, and we get
$$
\begin{aligned} \frac{1}{p} \int_\Omega b^\varepsilon|T_R(\omega^\varepsilon)(t)|^p
  + \varepsilon (p-1) \int_0^t ds \int_\Omega & b^\varepsilon |\nabla T_R(\omega^\varepsilon)|^2|\omega|^{p-2} \, dx
  \\
  & = \frac{1}{p} \int_\Omega b^\varepsilon|T_R(\omega_0)|^p \, dx.
\end{aligned}
$$
Therefore
$$\|b^{1/p} T_R(\omega)(t)\|_{L^p(\Omega)} \le \|b^{1/p} T_R(\omega_0)\|_{L^p(\Omega)}.$$
Letting $R$ go to $+\infty$, this gives
$$\|(b^\varepsilon)^{1/p} \omega^\varepsilon\|_{L^p(\Omega)}
   \le \|(b^\varepsilon)^{1/p} \omega_0\|_{L^p(\Omega)}.$$
and since $\omega_0 \in L^\infty$ and the estimate is uniform
$$\|\omega^\varepsilon\|_{L^\infty(\Omega)} \le \| \omega_0\|_{L^\infty(\Omega)}.$$
  Using the estimate obtained in the main estimate part we also get
$$ \|v^\varepsilon\|_{L^4(\Omega)} \le C(\|\omega^\varepsilon\|_{L^\infty(\Omega)} +
                                    \|\sqrt b v^\varepsilon\|_{L^2(\Omega)}).
                                     $$
                                     
  Next we note  that the set $\{b^\varepsilon\omega^\varepsilon, \varepsilon\in ]0, 1 ] \} $ is   relatively
compact   in $C^0([0,\infty); L_w^2(\Omega))$ and in
$C^0([0,\infty); L_{w*}^\infty(\Omega))$ as in \cite{LOT}.
  Thus $\{b^\varepsilon\omega^\varepsilon \}$ is a relatively compact set in $L^2_{loc}([0,\infty); H^{-1}(\Omega))$.
Therefore $\{\sqrt b v^\varepsilon\}$ is relatively compact in $L^2_{\rm loc}([0,\infty); L^2(\Omega))$.
This allows to pass to the limit in the viscous formulation and to get the global
existence of a weak solution of the inviscid stream-vorticity formulation.
\bigskip
\bigskip

\noindent{\bf Uniqueness of weak solutions.}

\begin{proof}
 Let $v_1$ and $v_2$ be two solutions of \eqref{Lake}.
 Then $v = v_1 - v_2$ satisfies
\begin{equation*}
\D_t v + v_2 \cdot \na v + \na p = - v \cdot \na v_1.
\end{equation*}
Therefore,
\begin{equation}
\label{estdif}
\frac{d}{dt} \big\Vert  \sqrt b v \big\Vert^2_{L^2}
\le 2 \int_{\Omega} b |v|^2 |\na v_1|
\end{equation}
  If $\na v_1  \in L^\infty$, this clearly implies that $v = 0$.
But in general, this $L^\infty$ estimate is not available, but
following \cite{Yo}, sharp $L^p$ estimates can be substituted.
  By Theorem \ref{mainest},
\begin{equation}
\label{b2}
C := \sup_{ p \ge 3} \  \Big\{
\frac{1}{p} \Big(\int_{\Omega}  | \na u_1 |^p \Big)^{\frac{1}{p}} \Big\} < \infty
\end{equation}
In the right hand side of \eqref{estdif}, we use Young's inequality
with
$$
  | \na v_1  | \in L^p ,  \quad
 b^{1/p'} |v |^{2/p'}  \in L^{p'}, \quad
  b ^{ 1 - 1/p' } | v |^{ 2 - 2/p'} =   | b v |^{  2/p} \in L^\infty,
$$
implying that
$y(t) :=  \big\Vert  \sqrt b v (t)  \big\Vert^2_{L^2} $ satisfies for all $p \ge 3$:
\begin{equation}
\label{ineq1}
\D_t y(t) \le p C M^{ 2 /p} \{ y(t)\}^{ 1/p'}
\end{equation}
with
$M = \Vert \sqrt b v  \Vert_{L^\infty} + \Vert \sqrt b v_2\Vert_{L^\infty}$.
Optimizing in $p$ yields
\begin{equation}
\label{ineq2}
\D_t y(t) \le   e C  y  (t)  \frac{1}{ \ln ( M^2 / y(t) )}.
\end{equation}
Since $y(0) = 0$, this implies that $y \equiv 0$.
\end{proof}

\bigskip
\bigskip
\bigskip
\noindent{\bf Acknowledgments.}
  The authors are partially supported by the french groupement de Recherche (GdR) "\'Equations d'Amplitudes et Propri\'et\'es Qualitatives" (EAPQ), managed by \'E. Lombardi, of the Centre National de Recherches Scientifiques (CNRS).
  D. {\sc Bresch} is also supported by a Rh\^one-Alpes fellowship obtained in 2004
on problems related to viscous shallow water equations and by the "ACI jeunes chercheurs 2004"
du minist\`ere de la Recherche "\'Etudes math\'ematiques de param\'etrisations en
oc\'eanographie". 


\end{document}